\documentclass{elsart}

\usepackage{amssymb, amsmath}

\newcommand{\al}{\alpha}
\newcommand{\be}{\beta}
\newcommand{\ga}{\gamma}
\newcommand{\de}{\delta}
\newcommand{\ze}{\zeta}
\newcommand{\io}{\iota}
\newcommand{\ka}{\kappa}
\newcommand{\la}{\lambda}
\renewcommand{\phi}{\varphi}
\newcommand{\si}{\sigma}

\newcommand{\Power}{\mathcal{P}}
\newcommand{\st}{\,|\,}
\newcommand{\forces}{\Vdash}
\newcommand{\sat}{\vDash}

\newcommand{\cat}{^{\,\frown}}

\renewcommand{\iff}{\leftrightarrow}

\newcommand{\restricted}{\upharpoonright}
\newcommand{\restr}{\!\upharpoonright\!}

\newcommand{\Add}{\textrm{Add}}

\DeclareMathOperator{\dom}{dom}

\newcommand{\GCH}{\textrm{GCH}}
\newcommand{\HOD}{\textrm{HOD}}

\newcommand{\ORD}{\textrm{Ord}}

\newcommand{\Succ}{\textrm{Succ}}
\DeclareMathOperator{\supp}{supp}
\newcommand{\trcl}{\textrm{trcl}}

\newcommand{\ZFC}{\textrm{ZFC}}

\newcommand{\ds}{\diamondsuit^*}
\newcommand{\dplus}{\diamondsuit^+}
\renewcommand{\diamond}{\diamondsuit}

\newcommand{\emin}[1]{\emph{#1}\index{#1}}

\begin{document}

\begin{frontmatter}

\title{Large Cardinals and Definable Well-Orderings of the Universe}

\author{Andrew D. Brooke-Taylor}

\address{Kurt G\"odel Research Center for Mathematical Logic,
University of Vienna, W\"ahringer Stra\ss e 25, 1090 Vienna, Austria}
\ead{andrewbt@logic.univie.ac.at}

\begin{abstract}
We use a reverse Easton forcing iteration to obtain a 
universe with a definable well-ordering,
while preserving the GCH and proper classes of a 
variety of very large cardinals.
This is achieved by coding using the principle $\ds_{\ka^+}$
at a proper class of cardinals $\ka$.  
By choosing the cardinals at which coding occurs sufficiently sparsely,
we are able to lift the embeddings witnessing the large cardinal 
properties without having to meet any non-trivial master conditions.
\end{abstract}

\end{frontmatter}

\nocite{*}

\section{Introduction}

A major theme of set theory in recent years has been
the construction of models of ZFC which contain various large cardinals while
at the same time enjoying properties analogous to those of 
G\"odel's constructible universe $L$. 
Generally this has been approached via the 
{\em inner model programme,} constructing canonical
$L$-like inner models for the 
large cardinals under consideration.  While this approach has had much 
success, there is a bound on the size of the large cardinals that have 
thus far been
accommodated by such techniques.  
For this
reason, Sy Friedman has proposed the {\em outer model programme,}
in which the goal is to construct $L$-like outer models containing 
large cardinals by the method of forcing.  
In doing so, one may obtain new consistency results for
large cardinals beyond the scope of current inner model theory.

One of the most striking properties of $L$ is the fact that it bears
a definable well-order, whose definition moreover involves no parameters.  
This may be expressed in terms of another well known inner model: 
$V$ has a well-order definable without parameters if and only if 
$V=\HOD$, the universe of
all hereditarily ordinal definable sets (in $V$).
It is well known that one may force to obtain such a model.
Indeed, McAloon \cite{McA:CROD} shows how to force a model of
$\ZFC+\GCH+V=\HOD+V\neq L$ starting from $L$, 
or using the continuum function for coding, a
model of $\ZFC+V=\HOD$ while preserving all measurables. 
Sy Friedman \cite{SDF:LCL} 
shows how to obtain a model of $\ZFC+\GCH+V=\HOD$ while preserving a 
(single, although the same argument works for boundedly many)
hyperstrong or $n$-superstrong cardinal, starting with techniques of Jensen to
make $V=L[A]$ for $A$ a subset of some sufficiently large cardinal.

In this article we improve upon these results, 
exhibiting a forcing construction that yields a universe satisfying both
$V=\HOD$ and the GCH, while preserving proper classes of a variety of large
cardinals.  The techniques will not in general suffice to preserve
all cardinals satisfying a given large cardinal property, 
but rather those satisfying
a combined large cardinal and anti-large cardinal axiom of a certain kind, 
such as ``superstrong but not a
limit of huge cardinals''.  
Of course, this is sufficient to prove the relative 
consistency of 
$$\ZFC+\GCH+V=\HOD+\exists\text{ a proper class of cardinals }\ka
\text{ such that }\phi(\ka)$$
for a variety of very strong large cardinal properties $\phi$.

The idea of our forcing is essentially to add unboundedly
large Cohen sets and 
use the fact that
every element of $V[G]$ is ordinal definable from
$A$, a class predicate for the added Cohen sets.  
Of course, we then
want to make $A$ itself definable in $V[G]$.  
To achieve this, we code up the
choices made by the generic in terms of whether or not 
some combinatorial principle holds at various cardinals.
Doing this while preserving the GCH puts a heavy constraint on 
which combinatorial principles can be used for such an encoding;
indeed, the GCH itself would otherwise be an ideal principle to use as a 
coding oracle, as in the work of McAloon \cite{McA:CROD}.
However, the existence of
$\diamondsuit^*_\ka$-sequences
also fits the bill nicely, without disturbing the GCH.

Coming at this from the other direction, we have a property
suitable to  be used as an oracle
(existence of $\ds_\ka$-sequences) and we want to use it to encode a
definable well-order of the (extension) universe.  
Instead of using some complicated
iteration with lots of bookkeeping, 
we may simply ``let the generic decide'' 
which way to force at each stage.  
This technique --- having an iteration at each stage of which the generic 
makes an initial decision that
determines the rest of the forcing poset at that stage ---
is not new in other
contexts;
see for example
Theorems
5.27 and 5.33 of \cite{SDF:FSCF}
and Section 3 of \cite{Ham:LP}.

Finally, we wish to do all this while preserving large cardinals.
As usual, we preserve the large cardinals in question by 
lifting the witnessing elementary embeddings.
Often this entails selecting a generic that lies below 
a specific \emph{master condition}, and so in general it is problematic 
to simultaneously preserve many large cardinals, as one will not in general
be able to choose a single
generic meeting all of the necessary master conditions.
The solution we present is to simply 
avoid the problem, making all of the master conditions in question trivial
by only coding at cardinals which will not lead to non-trivial
master condition requirements.  However, we may not be able to
achieve this for all large cardinals of the given kind, and this is what 
leads to the exceptions in our preservation theorems.

Throughout this article we shall assume that our ground model satisfies the
GCH.  It is known that that this can be forced while preserving a variety of
large cardinals; see for example \cite{SDF:LCL} which specifically
deals with the case of $n$-superstrong cardinals.

\section{The coding oracle $\ds_{\ka^+}$}\label{diamondStar}

Recall the following definition.

\begin{defn}
Let $\la$ be a regular cardinal, and let 
$D=\langle D_\al\st\al<\la\rangle$ be a sequence such that for
every $\al<\la$, $D_\al\subset\Power(\al)$ and $|D_\al|\leq|\al|$.
Then $D$ is said to be a 
\emph{$\ds_\la$-sequence}\index{$\ds_\la$!-sequence} if
for every $X\subset\la$, 
$\{\al\in\la\st X\cap\al\in D_\al\}$ contains a closed unbounded 
subset of $\la$.
The statement \emph{$\ds_\la$}\index{$\ds_\la$} 
is the statement that a $\ds_\la$-sequence exists.
\end{defn}

There are known $\ka^+$-closed, $\ka^{++}$-cc partial orders for forcing 
$\ds_{\ka^+}$ to hold or fail while preserving the GCH, 
for each infinite cardinal $\ka$. 
Indeed, in the $\ka=\omega$ case, the two directions
are given as exercises in \cite{Kun:ST}
(VII H.18--20 for forcing $\ds_{\omega_1}$ to hold, 
VIII J.3 for forcing it to fail).  
For completeness we outline the details (for the general case) here.

To force $\ds_\ka$ to hold there are at least two options.  Cummings,
Foreman and Magidor \cite{CFM:SSSR} present one possibility, with an iteration
in which the first iterand yields the eventual $\ds_{\ka^+}$-sequence
and the later iterands shoot clubs through $\ka^+$ to witness that the sequence 
is indeed a $\ds_{\ka^+}$-sequence.

For the more $L$-inclined, one can follow the lead of 
\cite{Kun:ST}, using a simpler forcing but a much more involved verification 
that it yields a $\ds_{\ka^+}$-sequence.  We use the following lemma.

\begin{lem}\label{DPEquivs}
Let $\ka^+$ be a successor cardinal.  
If there is a $\ka^+$-tree $T$ which is a subtree of ${}^{<\ka^+}2$, 
and an $h\in{}^{\ka^+}2$ such that
\begin{equation*}
\forall f\in{}^{\ka^+}2\,\exists g\in{}^{\ka^+}2
(f\in L(\{g,h\})\land 
g\text{ is the union of a cofinal branch in }T)
\end{equation*}
then $\dplus_{\ka^+}$ holds.
\end{lem}

Note that the conclusion here is the principle $\diamondsuit^+_{\ka^+}$,
which is stronger than $\ds_{\ka^+}$.  Indeed, we could code using 
$\diamondsuit^+_{\ka^+}$ rather than $\ds_{\ka^+}$, but have elected to 
use $\ds_{\ka^+}$ for simplicity.
The proof of Lemma \ref{DPEquivs} is similar to the proof of 
$\diamond^+_{\ka^+}$ in $L$; see \cite{Kun:ST} Exercise VI.9 for an outline
or \cite{Me:The} Section 4.1.1 for the full details.

With Lemma \ref{DPEquivs} at our disposal we may now simply force a 
Kurepa tree $T$ to exist and observe that it will satisfy the conditions
of the lemma.  We can do this with the following forcing, which is perhaps
slightly simpler than that presented (in the $\omega_1$ case) in \cite{Kun:ST}.

\begin{defn}
For any cardinal $\ka$, 
let $P_{\ka^+}$\index{$P_{\ka^+}$ ($\ds_{\ka^+}$ forcing)} 
denote 
the partial order whose elements are pairs $\langle X,\al\rangle$,
where $X\subset{}^{\ka^+}2$, $\al<\ka^+$, and $|X|\leq\ka$.
For $\langle X,\al\rangle,\langle Y,\be\rangle\in P_{\ka^+}$, 
say that $\langle Y,\be\rangle\leq \langle X,\al\rangle$
\index{$\leq$ relation! $\ds_{\ka^+}$ forcing} if and only if
$Y\supseteq X$, $\be\geq\al$, and for all $f\in Y$ there is a $g\in X$
such that $f\restr\al=g\restr\al$.
\end{defn}

The condition $\langle X,\al\rangle$
can be thought of as determining that the initial
segment up to level $\al$ of the ultimate $\ka^+$-Kurepa tree $T$ will 
consist of the $\al$-initial segments of the elements of $X$,
and further, forcing that
every element of $X$ will be the union of a cofinal branch of 
$T$.
Clearly $P_\ka$ is $\ka^+$-closed, has the $\ka^{++}$-cc, and has
cardinality $\ka^{++}$, so cardinals and the GCH are preserved.
To verify that the conditions of Lemma \ref{DPEquivs} hold for $T$ in the
generic extension $V[G]$ by $P_\ka$, we take $h\in({}^{\ka^+}2)^V[G]$
encoding $T$ itself and every element of $\Power(\ka)^{V[G]}$.
Every ground model function from $\ka$ to 2 appears as a cofinal branch of
$T$ up to bounded differences, every extension model function from 
$\ka$ to 2 has a nice name that can be encoded by a ground model function from
$\ka$ to 2, and $T$ determines the entire generic $G$, so indeed every
$f\in({}^{\ka^+}2)^{V[G]}$ lies in $L(\{g,h\})$ for some branch $g$ of $T$,
as required.

Making $\ds_{\ka^+}$ fail will simply be a matter of adding 
$\ka^{++}$-many Cohen subsets of $\ka^+$ and 
observing that $\ds_{\ka^+}$ will not hold in the extension.
We use the common notation $\Add(\la,\mu)$ for the partial order
with partial functions from $\la\times\mu$ to 2 of
cardinality less than $\la$ as conditions; thus, our forcing to make
$\ds_{\ka^+}$ fail will be $\Add(\ka^+,\ka^{++})$.
The arguments in this subsection actually work for any uncountable 
regular cardinal $\la$, so we present them at this level of generality, 
although for later sections we will only need the case when $\la=\ka^+$
for some $\ka$.

The partial order $\Add(\la,\la^+)$ is of course $\la$-closed, 
$\la^{+}$-cc, and has cardinality $\la^+$, so it preserves cardinals and the
GCH.
The verification that it 
destroys $\ds_\la$ 
is presented in \cite{Dev:VoD} for the case when $\la=\omega_1$,
but with appropriate modifications the argument can be lifted to work for 
any regular uncountable $\la$.  
We present this modified argument here.

To avoid the temptation to abuse notation, we make a definition to
recast $\ds_\la$-sequences in terms of characteristic functions,
giving notationally more convenient objects.

\begin{defn}\label{DSListing}
Let $\la$ be a regular cardinal and let 
$D=\langle D_\al|\al<\la\rangle$ be a sequence such that for every
$\al<\la$, $D_\al\subset\Power(\al)$ and $|D_\al|\leq|\al|$.
We say that $d$ is a \emin{listing} of $D$ if $d$ is a function 
on $\la$ with the property that for each $\al<\la$, $d(\al)$ enumerates the
characteristic functions of the elements of $D_\al$ in order type
$|D_\al|$.  That is,
\renewcommand\theenumi {\roman{enumi}}
\begin{enumerate}
\item for each $\al<\la$, $d(\al)$ is a a function from $|D_\al|$ to 
${}^\al2$; and
\item for all $\al<\la$ and $\be<|D_\al|$, there is some $S\in D_\al$
such that for all $\zeta<\al$, 
$d(\al)(\be)(\zeta)=1$ if and only if $\zeta\in S$; and
\item for all $\al<\la$ and $S\in D_\al$, there is a unique $\be<|D_\al|$
such that for all $\zeta<\al$, 
$d(\al)(\be)(\zeta)=1$ if and only if $\zeta\in S$.
\end{enumerate}
\renewcommand\theenumi {\arabic{enumi}}
\end{defn}

\begin{prop}
Suppose that $M\sat\ZFC+\GCH$ and $\la$ is a regular cardinal of $M$.  
Then if $G$ is $\Add(\la,\la^+)$-generic over
$M$, $M[G]\sat\ZFC+\GCH+\lnot\ds_{\la}$.
\end{prop}
\begin{pf}
We first claim that forcing with $\Add(\la,\la)$
destroys any $\ds_{\la}$-sequence of $M$.
Let $D$ be a $\ds_{\la}$-sequence of $M$.
For notational convenience, we may assume by expanding the sets
$D_\al$ if necessary that for every $\al<\la$, $|D_\al|=|\al|$.
Let $d$ be a listing of $D$; the assumption of the last sentence thus becomes
the statement that for each $\al<\la$, $d(\al)$ has domain $|\al|$.

Let $G_\la$ be $\Add(\la,\la)$-generic over $M$, let $\dot C$ be a 
name for a club subset of $\la$ in $M[G_\la]$,  
and let $\dot F$ name $\bigcup G_\la$.
We claim that in $M[G_\la]$, 
the subset of $\la$ with characteristic function $\bigcup G_\la$
is not correctly guessed by $D$ on all elements of $\dot C_{G_\la}$.
Suppose to the contrary that there is some $p\in G_\la$ such that
$$
p\forces
(\dot C\text{ is a club in }{\check \la})\land
\forall\al\in\dot C\exists\de<|\al|(\dot F\restricted\al=\check d(\al)(\de)).
$$
So that we may smoothly deal with the successor and inaccessible cases
simultaneously, let $\ga=\ka+1$ for $\ka$ such that $\ka^+=\la$ if $\la$
is a successor cardinal, and let $\ga=\la$ otherwise.
By induction 
on rank in ${}^{<\ga}2$, 
we may construct 
conditions 
\mbox{$\langle p_s\st s\in{}^{<\ga}2\rangle$}
and ordinals 
$\langle \al_s\st s\in{}^{<\ga}2\rangle$
such that the following properties hold.
\renewcommand\theenumi {\roman{enumi}}
\begin{enumerate}
\item\label{KillDSEmpty} 
$p_\emptyset\leq p$ and 
$p_\emptyset\forces\check\al_\emptyset\in\dot C$.
\item\label{KillDSOrd} 
For every $s\in{}^{<\ga}2$, $\dom(p_s)\in\la$.
\item\label{KillDSMono}
$s\subseteq t$ implies $p_t\leq p_s$ and $\al_t\geq\al_s$
\item\label{KillDSDiff}
For $i\in 2$ we have 
$\dom(p_{s\cat\langle i\rangle})\ni\dom(p_s)$, and
$$p_{s\cat\langle i\rangle}(\dom(p_s))=i.$$
\item\label{KillDSSucForc}
For $i\in 2$ we have $\al_{s\cat\langle i\rangle}>\dom(p_s)$ and
$p_{s\cat\langle i\rangle}\forces\check\al_{s\cat\langle i\rangle}\in\dot C$.
\item\label{KillDSJumpAls}
If $t\in{}^{\be}2$ then 
$\dom(p_t)\geq\sup(\{\al_s\st s\in{}^{\be}2\})$.
\item\label{KillDSLevels}
$\dom(s)=\dom(t)$ implies $\dom(p_s)=\dom(p_t)$.
\item\label{KillDSLims}
If $s\in{}^\eta2$ for $\eta$ is a limit ordinal, then 
$p_s=\bigcup_{\be<\eta}p_{s\restricted\be}$ and 
$$\al_s=\sup(\{\al_{s\restricted\be}\st\be<\eta\}).$$
\end{enumerate}
\renewcommand\theenumi {\arabic{enumi}}
Indeed, we may construct such 
\mbox{$\langle p_s\st s\in{}^{<\la}2\rangle$} and
$\langle \al_s\st s\in{}^{<\la}2\rangle$ by 
first extending $p$ as 
appropriate for (\ref{KillDSEmpty}) and (\ref{KillDSOrd}), at successor stages
extending to satisfy (\ref{KillDSDiff}), (\ref{KillDSSucForc}),
(\ref{KillDSJumpAls}) and (\ref{KillDSLevels}) in that order
while respecting (\ref{KillDSOrd}) and (\ref{KillDSMono}), and at limit 
stages satisfying (\ref{KillDSLims}).  Note that to satisfy 
(\ref{KillDSJumpAls}) at successor stages,
we rely on our assumption of the GCH to
give that $\{\al_s\st s\in{}^{\be}2\}$ is bounded below $\la$.

We claim that for all $s\in{}^{<\ga}2$, $p_s\forces\check\al_s\in\dot C$.
Of course from the definitions we need only check this for $s$ with
domain a limit ordinal.  But for such $s$,
$p_s\forces\check\al_{s\restricted\be}\in\dot C$ for all $\be<\dom(s)$,
and so since $p\geq p_s$ forces that $\dot C$ is a club and 
$\al_s=\sup(\{\al_{s\restricted\be}\st\be<\ga\})$,
$p_s\forces\check\al_s\in\dot C$.

Note that for any limit ordinal $\zeta<\ga$, $\al_s$ for $s\in{}^\zeta 2$
is independent of the choice of $s$: for $t$ with domain less than $\zeta$,
$$
\al_{t\cat\langle i\rangle}>\dom(p_{t})\geq
\sup(\{\al_s\st s\in{}^{\dom(t)}2\}),
$$
so
$$
\al_s=\sup(\{\al_{s\restricted\be}\st\be<\zeta\})
=\sup(\{\al_t\st\dom(t)<\zeta\}).
$$
Hence, let us denote $\al_s$ for $s\in{}^\zeta 2$ by $\al_\zeta$.
Observe further that because we have terms $\dom(p_t)$ 
interleaving with terms $\al_t$ in the above inequalities, and
$p_s=\bigcup_{\be<\zeta}p_{s\restricted\be}$ for 
$s\in{}^\zeta 2$, we have $\dom(p_s)=\al_\zeta$ for $s\in{}^\zeta2$.

But now let $\mu$ be the least cardinal such that $|\al_\mu|=\mu$; 
such a $\mu<\la$ can easily be found by a typical closure argument.
For each $s\in{}^{\mu}2$ we have a condition $p_s$ such that
$p_s\restricted\la\in{}^{\al_\mu}2$.  Moreover, (\ref{KillDSDiff}) dictates
that for $s\neq t\in{}^{\al}2$, $p_s\restricted\la\neq p_t\restricted\la$.
Thus, we have $2^\mu$ distinct elements of ${}^{\al_\mu}2$, so not all
of them can be of the form $d(\al_\mu)(\de)$ for $\de<|\al_\mu|=\mu$.
So let $s\in{}^{\mu}2$ be such that for all $\de<\mu$,
$p_s\restricted\la=p_s\restricted\al_\mu\neq d(\al_\mu)(\de)$.
But then 
$$
p_s\forces(\check\al_\mu\in\dot C)\land
\forall\de<|\check\al_\mu|
(\dot F\restricted\al_\mu=\check{(p_s\restricted\al_\mu)}
\neq\check d(\check\al_\mu)(\de))
$$ 
contradicting the fact that $p_s\leq p$.
We have therefore shown that $\Add(\la,\la)$ 
destroys any ground model $\ds_\la$-sequences.

Next we claim that all ground model $\ds_\la$-sequences will be destroyed
by our forcing $\Add(\la,\la^+)$.  Since 
$\Add(\la,\la^+)\cong\Add(\la,\la)*\Add(\la,\la^+)$, 
it suffices to show that such a $D$ cannot be resurrected 
after being killed by the initial 
$\Add(\la,\la)$ piece of the forcing.  Now the statement that $D$ is
not a $\ds_\la$ sequence is equivalent to there being a subset $X$ of $\la$
such that the set of $\al$ such that $X\cap\al\notin D_\al$ is 
stationary in $\la$.  The partial order 
$\Add(\la,\la^+)$ preserves stationary subsets of
$\la$, so $D$ will continue not being a $\ds_\la$-sequence after 
subsequently forcing with $\Add(\la,\la^+)$.  Hence, $D$ is not 
resurrected, and we may conclude that all ground model $\ds_\la$-sequences
are destroyed by the forcing $\Add(\la,\la^+)$.

Now suppose that some $D$ of the right form to be a $\ds_\la$-sequence
(that is, satisfying the assumptions of Definition \ref{DSListing}) is
added by $\Add(\la,\la^+)$; we wish to show that $D$ is not in fact a 
$\ds_\la$-sequence in $M[G]$.  
Let $d\in M[G]$ be a listing of $D$.
Since $\Add(\la,\la^+)$ is $\la$-closed and hence adds no new
$<\la$-tuples of ground model sets, $d(\al)$ is an element of $M$ for each
$\al<\la$.  Therefore, $d$ can be named by a name $\dot d$ which involves
for each $\al<\la$
a single antichain of $\Add(\la,\la^+)$ to determine $d(\al)$.
To be precise, if $A_\al$ is a maximal antichain of conditions that
determine $d(\al)$, and for $p\in A_\al$ 
we write $f_p$ for that function such that
$p\forces\dot d(\check\al)=\check f_p$, then we may take $\dot d$ to be
$$
\dot d=\bigcup_{\al<\la}
\Big\{\Big\langle\big\{\langle\check\al,1\rangle\big\},1\Big\rangle,
\Big\langle\big(\big\{\langle\check\al,1\rangle\big\}\cup
\big\{\langle\check f_p,p\rangle\st p\in A_\al\big\}\big),1\Big\rangle
\Big\}.
$$
Since $\Add(\la,\la^+)$ has the $\la^+$-chain condition,
$|\bigcup_{\al<\la}A_\al|\leq\la$, and so
$$
\Big|\bigcup_{\al<\la}\bigcup_{p\in A_\al}\dom(p)\Big|\leq\la.
$$
Thus, there is some common upper bound $\ga<\la$ on the domains of
those conditions $p$ appearing in $\dot d$.
Now 
$$\Add(\la,\la^+)\cong\Add(\la,\ga)\times\Add(\la,\la^+),$$
and if $G_\ga=G\cap\Add(\la,\ga)$, we have
$d\in M[G_\ga]$.  
Since $\Add(\la,\ga)$ is $\la$-closed, 
$$\Add(\la,\la^+)^{M[G_\ga]}=\Add(\la,\la^+)^M.$$
So by what we have already shown, 
if $G^\ga=G\cap\Add(\la,\la^+)$, 
then $d$ does not represent a 
$\ds_\la$-sequence in $M[G_\ga][G^\ga]=M[G]$.
Therefore, there are no $\ds_\la$-sequences in $M[G]$.
\hfill\qed\end{pf}

\section{Forcing a definable well-order}\label{ForDWO}

In this section we exhibit our forcing which yields a universe with a 
definable well order.
There is much flexibility in the definition we shall present,
a fact which we will later exploit when trying to preserve various different
kinds of large cardinals.

As discussed in the introduction, the general idea of our forcing 
is to use 
$\ds_\ka$ at various $\ka$ to act as an oracle,
coding up a proper class of
ordinals from which our well-order will be defined.
In fact, we further obtain that the extension $V[G]$ is of the
form $L[A]$ for $A$ a definable class in $V[G]$.  In some sense this
is the closest to $L$ we can hope to get while trying to preserve very
large cardinals --- it follows from Kunen's theorem that 
$V\neq L(x)$ for any set $x$ if $V$ contains strong cardinals, 
and of course $A$ cannot be taken to be definable over $L$ as that would
give $L[A]=L$.  On the other hand, it is possible to have properties
very different from those of $L$ coded into $A$
(for example, the failure of $\ds_{\ka^+}$ for many cardinals $\ka$!), 
so in itself this should not
be thought of as a resolution of the outer model programme.

We wish to force at various successor cardinals to 
``switch $\ds$ on or off'', 
and then use this as an oracle to make the universe well-orderable.  
Perhaps the most natural sequence of cardinals at which to do this would
be simply the class of \emph{all} infinite successor cardinals.  
However, for consideration of large cardinal preservation, it will be 
convenient to use more restricted classes of successor cardinals,
and we present our results in this generality.

\begin{defn}
A definable class $C$ of cardinals is a \emph{coding class} if
there is a definable class $B$ of cardinals such that
\begin{enumerate}
\item $C$ is a proper class, and
\item every element of $C$ is a successor cardinal, and
\item if $B$ is a set, then every successor cardinal greater than the 
supremum of $B$ is in $C$, and
\item for every element $\be$ of $B$, $\be^+\in C$, and
the least successor cardinal greater than $\be$ that is not in $C$ 
(if such exists) is the successor of an inaccessible cardinal, and
\item $B$ is countably closed.
\end{enumerate}
\end{defn}
Given a coding class $C$, we will denote by $c$ the increasing enumeration
of $C$.
For our present purposes, one may think of $C$ being the class of
\emph{all} successor cardinals, with $B$ empty and 
$c$ being the function
$\aleph_{\,\cdot\,+1}:\al\mapsto\aleph_{\al+1}$.  
Of course there are concerns regarding the
absoluteness of $C$ which we shall address in due course; 
unless otherwise stated, $c(\al)$ should be taken to 
be computed in the ground model $V$.

So let $C$ be a fixed coding class, with $c$ its increasing enumeration.
We retain the notation of Section \ref{diamondStar} of
$P_{\ka^+}$ being our forcing to produce a $\ds_{\ka^+}$-sequence.  
For ease of notation let us set 
$Q_{\ka^+}=\Add(\ka^{+},\ka^{++})$, the forcing that quashes all 
$\ds_{\ka^+}$-sequences.  Note again that we assume $V\sat\GCH$. 
For each ordinal $\al$, 
let $R_{c(\al)}$ be the sum of $P_{c(\al)}$ and $Q_{c(\al)}$, that is,
the partial order given by
combining disjoint copies of $P_{c(\al)}$ and $Q_{c(\al)}$ below a new 
maximum element in the obvious way.
For concreteness, let us set ${1}_{R_{c(\al)}}=\emptyset$, and let
$$
R_{c(\al)}=\{{1}_{R_{c(\al)}}\}\cup
(\{0\}\times P_{c(\al)})\cup(\{1\}\times Q_{c(\al)}).
$$
\index{$R_\ka$ (iterand of $\ds$ Oracle Forcing)}
For $r_0, r_1\in R_{c(\al)}$, $r_1\leq r_0$ 
\index{$\leq$ relation! $\ds$ Oracle Forcing iterand}
if and only if either 
$r_0={1}_{R_{c(\al)}}$, 
or $r_0=\langle i,r'_0\rangle$ and $r_1=\langle i,r'_1\rangle$ 
for some $i\in 2$ and
$r'_1, r'_0\in P_{c(\al)}\cup Q_{c(\al)}$ such that 
$r'_1\leq_{P_{c(\al)}}r'_0$ or 
$r'_1\leq_{Q_{c(\al)}}r'_0$.
Clearly $R_{c(\al)}$ will have cardinality $c(\al)^+$, be $c(\al)$-closed, 
and have the $c(\al)^+$-cc, 
since these statements are true of both $P_{c(\al)}$ and $Q_{c(\al)}$.
For $\ga\in\ORD$ not of the form $c(\al)$ for some $\al$,
let $R_\ga$ be the trivial forcing.  

\begin{defn}\label{DSOracleForc}
The \emin{$\diamond^*$ Oracle Partial Order} $S$ 
\index{$S$ ($\ds$ Oracle forcing)} is 
the reverse Easton iteration of $\dot R_\al$ as above for $\al\in\ORD$.
\index{partial order!$\diamond^*$ Oracle}
\end{defn}
Note that with only trivial forcings used between cardinals, 
Easton support is the same as taking direct limits at inaccessibles
and inverse limits elsewhere.
In particular, this implies that $S$ can be factored as
$S_\al*\dot S^\al$ for any stage $\al$ --- see for example
Lemma 21.8 of \cite{Jech:ST}.

\begin{lem}\label{DWOForcPres}
If $V\sat\ZFC+\GCH$ and $G$ is generic for the 
$\ds$ Oracle Partial Order $S$ over $V$, then
$V[G]$ satisfies $\ZFC+\GCH$ and has the same cardinals as $V$.
\end{lem}
\begin{pf}
As is generally the case for reverse Easton iterations used in practice,
$S$ is tame because the iterands are increasingly closed
(see Lemmata 2.22 and 2.31 of \cite{SDF:FSCF}) so ZFC is preserved.
To prove that cardinals and the GCH are preserved, we argue by induction
on the length of the iteration.
Successor stages are immediate from the fact that $R_{\ka^+}$ is
$\ka^+$-closed and \mbox{$\ka^{++}$-cc.}
For limit stages $\la$, cardinals and the GCH are
preserved below $\la$ 
by the closure of the tail parts of the iteration $S_\la$.
There is a dense suborder of $S_\la$ of size at most $\la^+$ if $\la$ is 
singular or $\la$ is $\la$ is regular, so the GCH and cardinals are
preserved above
$\la^+$, or $\la$ in the regular case.  
It therefore only remains to show that in the $\la$ singular case,
the GCH holds at $\la$ and $\la^+$ is preserved, and this follows by
considering nice names for subsets of $\la$ built up as the
union of nice names for subsets of smaller cardinals.
\hfill\qed\end{pf}

In particular, note that the class of successor cardinals is unchanged at
each stage of the iteration, so if $c=\aleph_{\cdot+1}$, then $c$ is
absolute.

Considering the factorisation of $S$ as $S_\ka*S^\ka$, with
$S^\ka$ being $\ka$-closed, also gives the following.

\begin{lem}\label{DWOForcInacc}
Forcing with $S$ preserves inaccessible cardinals.\hfill\qed
\end{lem}

Another basic property of this forcing we shall need is that
after applying it, $\ds$ holds at exactly those points in the sequence $c$
where we expect it to: those $c(\al)$ such that
$\langle0,1_{P_{c(\al)}}\rangle\in G(c(\al))$.

\begin{lem}\label{DSWhereExpected}
Let $V\sat\ZFC+\GCH$ and let $G$ be $S$-generic over $V$.  Then
for every $\al\in\ORD$, $V[G]\sat\ds_{c(\al)}$ if and only if
$\langle0,1_{P_{c(\al)}}\rangle\in G(c(\al))$.
\end{lem}
\begin{pf}
Let $\ka^+$ be of the form $c(\al)$, and consider the factorisation of
$S$ as $S_{\ka^+}*R_{\ka^+}*S^{\ka^++1}$.  
Clearly $V[G_{\ka^+}*G(\ka^+)]\sat\ds_{\ka^+}$ if and only if
$\langle0,1_{P_{c(\al)}}\rangle\in G(c(\al))$,
and $S^{\ka^++1}$ is $\ka^{++}$-closed,
so any $\ds_{\ka^+}$-sequence of $V[G_{\ka^+}*G(\ka^+)]$ remains a 
$\ds_{\ka^+}$-sequence of $V[G]$, 
and no new $\ds_{\ka^+}$-sequences are added by $G^{\ka^++1}$.
\hfill\qed\end{pf}

\begin{thm}\label{DWOForcWorks}
Let $V\sat\ZFC+\GCH$, and let $S$ be the $\ds$ Oracle Partial Order
as defined above.
If $G$ is $S$-generic over $V$, then 
there is a definable class of ordinals $A$ of $V[G]$ such that 
$V[G]=L[A]$.
In particular, 
$V[G]=\HOD^{V[G]}$, and  
there is a definable well-order on $V[G]$.
\end{thm}
\begin{pf}
The class $A$ will of course be 
$\{\al\in\ORD\st\ds_{c(\al)}\text{ holds}\}.$
Clearly $A$ is definable in $V[G]$; the question will be what the
relationships between $A$, $C^{V[G]}$ and $C^V$ are, because of course
$C$ need not be absolute.

Actually in many specific cases of interest, $C$ \emph{will} be absolute.
As mentioned above, the class of all successor cardinals will be absolute,
and in context of the next section one could get absoluteness from
the fact that we are preserving the large cardinals from which $B$ is
defined, along with Hamkins' Gap Forcing Theorem \cite{Ham:Gap}.
But in any case, the requirements we have placed on $B$ will 
give $C^V$ and $C^{V[G]}$ sufficient agreement
to show that every set in $V[G]$ is encoded into $A$.
Namely,
since $V[G]$ satisfies the Axiom of Replacement with respect to 
formulas involving a predicate for $V$
(see \cite{SDF:FSCF} Lemma 2.19), 
the usual argument to show that countably
closed unbounded sets have 
countably closed unbounded intersection goes through for $B^V\cap B^{V[G]}$, 
and we see that there are unboundedly many 
cardinals in $B^V\cap B^{V[G]}$.  
We claim that this agreement is
sufficient for our purposes.

So suppose $x\in V[G]$, and let $X$ be a subset of $\mu=|\trcl(\{x\})|$ 
coding up $x$.
We claim that for every 
$\be$ such that $c(\be)\geq\mu$,
$X$ appears in the choices made by the generic between $c(\be)$ and the least 
inaccessible greater than $c(\be)$ (if one exists),
where in each case $c(\be)$ is to be computed in $V$.  
To see this, 
let $\iota_{c(\be)}$ denote the least inaccessible greater than $c(\be)$
if such exists or $\ORD$ otherwise, and consider the class
\begin{eqnarray*}
D_{X,\be}&=&\big\{s\in S^{\mu+1}\st\exists\ga\in\iota_{c(\be)}\,
\big(\ga\geq\be\land\forall\al<\mu\\
&&\,\ \qquad\qquad\big(\al\in X\rightarrow 
\ \forces_{S^{\mu+1}}s(c(\ga+\al))\leq_{\dot R_{c(\ga+\al)}}
\langle 0,\dot 1_{P_{c(\ga+\al)}}\rangle \land\\
&&\,\ \qquad\qquad\,\ \al\notin X\rightarrow 
\ \forces_{S^{\mu+1}}s(c(\ga+\al))\leq_{\dot R_{c(\ga+\al)}}
\langle 1,\dot 1_{Q_{c(\ga+\al)}}\rangle\ \big)\big)\big\}.
\end{eqnarray*}
Because a direct limit is taken at $\iota_{c(\be)}$, 
we have for any $s\in S^{\mu+1}$
that \mbox{$\supp(s)\cap\iota_{c(\be)}$} 
is bounded in $\iota_{c(\be)}$, so we may 
extend $s\restr\iota_{c(\be)}$ to an element of 
$D_{X,\be}\cap S^{[\mu+1,\iota_{c(\be)})}$, and then 
``re-attach the tail of $s$'' to get an extension of $s$ in $D_{X,\be}$.
Hence, for each $\be$ with ${c(\be)}\geq\mu$, 
the class $D_{X,\be}$ is dense in $S^{\mu+1}$,
and so has non-empty intersection with $G^{\mu+1}$.

Now because inaccessibles are absolute between $V$ and $V[G]$ 
(Lemma~\ref{DWOForcInacc}), 
if $\la\in B^V\cap B^{V[G]}$ is greater than $\mu$ and
$\io_\la$ is the least inaccessible greater than $\la$, then 
$[\la,\io_\la)\cap\Succ\subset C^V\cap C^{V[G]}$, where $\Succ$
denotes the class of successor cardinals.
Therefore, taking $\be$, $s$ and $\ga$ such that 
$\be$ is least with $c(\be)\geq\la$,
$s\in D_{X,\be}\cap G^{\mu+1}$,
and $\ga$ is as in the definition of $D_{X,\be}$ 
witnessing that $s\in D_{X,\be}$,
we have that $c^V``[\ga,\ga+\mu)=c^{V[G]}``[\ze,\ze+\mu)$ 
for some ordinal $\ze$, and so indeed, $X\in L[A]$.
Hence, we have shown that $V[G]=L[A]$, as required.
\hfill\qed\end{pf}

\section{Preserving large cardinals}

As mentioned in the introduction, if we wish to preserve large cardinals
while forcing, we will generally have master conditions to be hit by our
generic, which will be problematic if we wish to preserve many large cardinals.
In other settings this can be overcome by making the partial order
sufficiently homogeneous that generics containing particular master conditions
can be constructed in \emph{any} generic extension --- see for example
\cite{Me:LCMor} for the case of forcing gap-1 morasses to exist at every 
regular cardinal.  
We need another approach, however, as
our forcing partial order is inherently inhomogeneous --- indeed,
any forcing that yields a model of $V=\HOD$ must be inhomogeneous,
as $\HOD$ of the generic extension by a homogeneous forcing must 
be contained in $V$
(see \cite{Kun:ST}, Exercise VII E1).

The solution to this problem comes from the extra flexibility we have
because we are interested in forcing a global principle, rather than
a local principle at, say, every regular cardinal.  We can ``thin out'' our 
forcing partial order, still obtaining a definable well-order of the
extension universe,
but finessing the issue of master conditions
by making the forcing trivial at every point where master conditions 
might be required.

To facilitate this thinning out, we make the following definitions.
\begin{defn}\label{LCBd}
Suppose that $\phi$ is a formula in one variable, and more specifically, that:
\begin{enumerate}
\item\label{phiMeas} $\phi(\ka)\equiv$ ``$\ka$ is measurable'', or
\item for some ordinal $\eta$, $\phi(\ka)\equiv$ ``$\ka$ is $\eta$-strong'', or
\item $\phi(\ka)\equiv$ ``$\ka$ is Woodin'', or
\item for some $n\in\omega+1$, $\phi(\ka)\equiv$ ``$\ka$ is $n$-superstrong'', or
\item $\phi(\ka)\equiv$ ``$\ka$ is hyperstrong'', or
\item for some definable function $g$, 
$\phi(\ka)\equiv$ ``$\ka$ is $g(\ka)$-supercompact'', or
\item for some ordinal $\eta$, $\phi(\ka)\equiv$ ``$\ka$ is $\eta$-extendible'', or 
\item\label{phimHuge} for some $m\in\omega$, $\phi(\ka)\equiv$ ``$\ka$ is $m$-huge''.
\end{enumerate}
A cardinal $\la$ is a \emph{$\phi$-bound}
if
$\la$ is an infinite cardinal,
$\la$ is not Mahlo, and
if $\phi(\ka)$ holds for any $\ka<\la$,
then there is an elementary embedding 
$j$ with critical point $\ka$ witnessing the fact that $\phi(\ka)$ holds, 
such that 
\renewcommand\theenumi {\roman{enumi}}
\begin{enumerate}
\item if $\phi(\ka)\equiv$ ``$\ka$ is measurable'' then $\ka<\la$ 
(that is, no extra requirement),
\item if $\phi(\ka)\equiv$ ``$\ka$ is $\eta$-strong''
 then $\beth_{\ka+\eta}^+<\la$,
\item if $\phi(\ka)\equiv$ ``$\ka$ is Woodin''
 then for all $f:\ka\to\ka$ there is an $\al\in\ka$ and a $j:V\to M$
elementary such that $f``\al\subset\al$, $\textrm{crit}(j)=\al$, 
$V_{j(f)(\al)}\subseteq M$, and $\beth_{j(f)(\al)}^+<\la$.
\item if $\phi(\ka)\equiv$ ``$\ka$ is $n$-superstrong''
 then $\beth_{j^{n}(\ka)}<\la$,
\item if $\phi(\ka)\equiv$ ``$\ka$ is hyperstrong''
 then $\beth_{j(\ka)+1}<\la$,
\item if $\phi(\ka)\equiv$ ``$\ka$ is $g(\ka)$-supercompact''
 then $ g(\ka)^{<\ka}<\la$,
\item if $\phi(\ka)\equiv$ ``$\ka$ is $\eta$-extendible''
 then $ \zeta<\la$ for the 
$\zeta$ such that $j:V_{\ka+\eta}\to V_\zeta$, and
\item if $\phi(\ka)\equiv$ ``$\ka$ is $m$-huge''
 then $ 2^{j^m(\ka)}<\la$.
\end{enumerate}
We say that $\la$ is a 
\emin{minimal $\phi$-bound}\index{$\phi$-bound!minimal}
if for every cardinal $\nu<\la$ there is a cardinal $\ka$ which
is not a $\phi$-bound such that $\nu\leq\ka<\la$.
\end{defn}

Clearly one may construct minimal $\phi$-bounds by the usual
methods for obtaining fixed points.  
The non-Mahloness assumption is simply a convenient way to ensure that
the large cardinals were are interested in are not considered 
$\phi$-bounds.  Also note that we have stuck with $\beth$ notation
for clarity in the arguments to come,
even though we are assuming the GCH and so $\aleph$ notation would be 
equivalent.

Unless otherwise specified, 
$\phi$ shall henceforth denote one of the listed large cardinal properties
\ref{phiMeas}--\ref{phimHuge}, and
for convenience we shall refer to cardinals $\ka$ satisfying $\phi(\ka)$
as \emph{$\phi$-cardinals}\index{$\phi$-cardinal}.
This list of large cardinal properties, which will be the 
ones that are preserved in Theorem~\ref{DWOClassLCs},
should by no means be thought of as encompassing all 
large cardinals for which the
techniques of this chapter are applicable.
Rather, it is a representative list of well-known large cardinals 
each witnessed by boundedly many elementary embeddings
so that bounds may be constructed for them.
For reasons that will become apparent, it would also be of interest 
(and moreover straightforward) to include
large cardinals of the form ``$\phi$ a limit of $\phi$'' 
(for example, a measurable limit of measurables).  
However, if we do not wish to move to a more general statement,
we must draw the line somewhere!

Once we have a
$\phi$-bound, the succeeding cardinals will remain $\phi$-bounds
for some time.  
The following lemma in this direction will be sufficient for our 
purposes.
\begin{lem}\label{InaccLLCB}
Suppose $\aleph_\be$ is a $\phi$-bound.
Then for all $\ga$ less than the least inaccessible greater than
$\aleph_\be$ (or all $\ga$ if no such inaccessible exists), 
$\aleph_{\be+\ga}$ is a $\phi$-bound.
\end{lem}
\begin{pf}
This is immediate from the fact that the least cardinal $\ka$
which is \emph{not}
a $\phi$-bound above a given $\phi$-bound will satisfy $\phi(\ka)$,
and hence be inaccessible.
\hfill\qed\end{pf}

With this fact, we are ready to define the points at which we shall
perform our coding.
\begin{defn}\label{CodingPt}
A cardinal $\la$ is a \emph{$\phi$-coding point}\index{coding point, $\phi$-} 
if 
\begin{enumerate}
\item $\la$ is a successor cardinal, and 
\item $\la$ is a $\phi$-bound, and
\item if there is a cardinal $\ka>\la$ such that $\phi(\ka)$, 
there is a minimal $\phi$-bound $\mu\leq\la$ such that
$\la$ is less than the least inaccessible cardinal greater than $\mu$.
\end{enumerate}
\renewcommand\theenumi {\roman{enumi}}
\end{defn}
Thus, our coding points come after each minimal $\phi$-bound,
going on until the next inaccessible cardinal, or indefinitely if there is
no next $\phi$-cardinal.
Clearly the class $C$ of all $\phi$-coding points is a coding class,
witnessed by the class $B$ of minimal $\phi$-bounds.

\begin{thm}\label{DWOClassLCs}
Let $V\sat\ZFC+\GCH$, and let $\phi$ be one of the large cardinal
properties \ref{phiMeas}--\ref{phimHuge} listed in Definition~\ref{LCBd}.
Let $C$ be the coding class of all $\phi$-coding points of $V$,
and let $S$ be the $\ds$ Oracle Partial Order defined from $C$.
Suppose $G$ is $S$-generic over $V$.  Then there is 
$V[G]$-definable class of ordinals $A$ such that $V[G]=L[A]$.
Further, if $\ka$ is a $\phi$-cardinal in $V$ that is not a limit of
$\phi$-cardinals, then in $V[G]$ $\ka$ remains a $\phi$-cardinal.
\end{thm}
\noindent So for example, we can preserve all measurable cardinals that are not
limits of measurables; 
see below for a discussion of extensions 
strengthening this.
\begin{pf}
As before, we denote by $c$ the increasing enumeration of $C$, and let
$B$ denote the class of minimal $\phi$-bounds in $V$.
Theorem \ref{DWOForcWorks}
gives that $V[G]=L[A]$, and it only
remains to show that any $\phi$-cardinal $\ka$ of $V$ that
is not a limit of $\phi$-cardinals remains a $\phi$-cardinal in $V[G]$.
We prove this by lifting embeddings witnessing
$\phi(\ka)$, and moreover taking these embeddings to be
given by ultrapowers or extenders where appropriate.
This will allow us to use 
representation results about the codomain of such embeddings;
a good reference for these is \cite{Kan:THI}.
We deal with each large cardinal property separately.\\

\noindent{\bf Measurable Cardinals.}
Let $j:V\to M$ be an ultrapower embedding witnessing the measurability of $\ka$
with $j(\ka)$ least.  
We shall 
construct an $S^M$-generic $G^*$ over $M$ in $V[G]$, such that we can
lift $j$ to $j^*:V[G]\to M[G^*]$.
Note that the $\phi$-coding points of $V$ less than $\ka$ 
are in fact bounded below $\ka$ since the class of 
measurable cardinals is.
Hence, $S^M$ is $S^V$ up to stage $\ka$, and is trivial from $\ka$ to $j(\ka)$.
We may therefore take $G^*_{\ka}=G_{\ka}$, trivially extend to $G^*_{j(\ka)}$,
and have a lift of $j$
to $j':V[G_{\ka}]\to M[G^*_{j(\ka)}]$.
To define $G^{*j(\ka)}$, 
note that every element of $M$
has the form $j(f)(\ka)$, 
where $f:\ka\to V$ is a function in $V$, and so every element of
$M[G^*_{j(\ka)}]$ has the form $\si_{G_{j(\ka)}}$, where 
$\si$ has the form $j(f)(\ka)$.
We claim that the filter on $S^{j(\ka)}$ generated by
$j'``G^{j(\ka)}$ is $(S^{j(\ka)})^M$-generic over $M$.

So suppose that $D$ is a dense class in $M[G^*_{j(\ka)}]$, 
defined (in $M[G^*_{j(\ka)}]$) relative to the parameter $d\in M[G^*_{j(\ka)}]$
by $D=\{x\st\psi(x,d)\}$.
Let $\sigma$ be an $S_{j(\ka)}^M$-name in $M$ such that
$d=\si_{G^*_{j(\ka)}}$, and let $f:\ka\to V$ in $V$ be such that
$\si=j(f)(\ka)$.
Since $\ka$ is not a measurable-coding point,
$(S^{\ka})^V$ is $\ka^+$-closed, and we see that
it is dense for $s\in(S^{\ka})^V$ to
extend an element of 
the class $D_\al=\{x\st\psi(x,f(\al)_{G_\ka})\}$ of $V$
whenever $\al\in\ka$ with $f(\al)$ an $S_\ka$-name and
$D_\al$ dense in $(S^{\ka})^V$.
Therefore, we may take such
an $s$ lying in $G^{\ka}$.  By elementarity, it follows that
$j(s)$ extends an element of $D$. Hence, the filter generated by 
$j'``G^{\ka}$ is indeed $(S^{j(\ka)})^{M[G^*_{j(\ka)}]}$-generic over $M$.
By the Lifting Lemma, it follows that there is an elementary embedding
$j^*:V[G]\to M[G^*]$ lifting $j$, and so $\ka$ is measurable in $V[G]$.\\

\noindent{\bf $\eta$-Strong Cardinals.} 
We may assume that our $\eta$-strong embedding $j:V\to M$ is
an extender ultrapower embedding, with every element of $M$ having the
form $j(f)(a)$, with $a$ a finite tuple from $|V_{\ka+\eta}|^+$ and
$f$ a function in $V$ from $[\ka]^{|a|}$ to $V$.
As in the measurable cardinal case, we show that the filter generated by 
$j'``G^\ka$ is $(S^{j(\ka)})^{M[G^*_{j(\ka)}]}$-generic,
observing that
by the definition of a $\phi$-bound, $S^{\ka}$ is trivial
up to at least stage $|V_{\ka+\eta}|^{++}$, 
and so we have the requisite 
closure to make the argument go through. 
We can therefore lift $j$ to $j^*:V[G]\to M[G^*]$. 
We may also conclude that $V^{V[G]}_{\ka+\eta}\subseteq M[G^*]$
from a nice names argument, since $V^V_{\ka+\eta}\subseteq M$, 
$G_{\ka+\eta}=G^*_{\ka+\eta}$, 
and $S_{\ka+\eta}$ is trivial beyond some bound below 
$\ka$.  Hence, $\ka$ is $\eta$-strong in $V[G]$.\\

\noindent{\bf Woodin Cardinals.}
The situation for Woodin cardinals is somewhat different from 
that for the 
other large cardinals listed here, since Woodinness is witnessed by
multiple embeddings.  However, since there are only boundedly many, 
this point will not present a problem. 
Let $\iota$ denote the supremum of the Woodin-coding points less than
$\ka$.  Note that for any two functions $f,g:\ka\to\ka$,
if $f(\al)\leq g(\al)$ for all $\al\in\ka$, then
$$
\{\al\in\ka\st f``\al\subseteq\al\}\subseteq
\{\al\in\ka\st g``\al\subseteq\al\}.
$$
Now given a name $\dot f$ for a function from $\ka$ to 
$(\iota,\ka)$, we can find a function $\bar f:\ka\to(\iota,\ka)$ in $V$
such that $1_S\forces\dot f\leq\check{\bar f}$, since the forcing 
iterands are trivial from $\iota$ to $\ka$.
To prove that Woodinness is preserved, then, 
it is sufficient to show that for every function
$f$ in $V$ from $\ka$ to the interval $(\iota,\ka)$, 
there is an $\al$ with $f``\al\subseteq\al$
and a $j^*$ from $V[G]$ to $N$ an inner model of $V[G]$ such that
$\textrm{crit}(j^*)=\al$ and $V_{j^*(f)(\al)}\subseteq N$.
Since $\ka$ is Woodin in $V$, we have for each such $f$ an $\al$ with
$f``\al\subseteq\al$ and an elementary embedding $j:V\to M$ such that
$\textrm{crit}(j)=\al$ and $V_{j(f)(\al)}\subseteq$; that is, an
$\eta$-strong embedding for $\eta$ such that $\al+\eta=j(f)(\al)$.
But this is simply a case of $\eta$-strength, so we can lift $j$ 
to $j^*$ as above.
Since $j^*\restr V=j$, $j^*(f)(\al)=j(f)(\al)$, 
and we are done.

\noindent{\bf $n$-Superstrong Cardinals and Hyperstrong Cardinals.}
The argument is analogous to the measurable and $\eta$-strong cases, 
this time with our extender models 
having elements of the form $j(f)(a)$ with $a$ in
$V_{j^n(\ka)}$ and $f$ with domain $V_{j^{n-1}(\ka)}$ 
in the case of $n$-superstrong
cardinals, and $a\in V_{j(\ka)+1}$ and $f$ with domain $V_{\ka+1}$
in the case of hyperstrong cardinals.
The required level of agreement between $V[G]$ and $M[G^*]$ 
again follows from a nice names argument, noting in the
hyperstrong case that $V_{j(\ka)+1}\in M\iff H_{j(\ka)^+}\in M$.\\

\noindent{\bf $g(\ka)$-Supercompact Cardinals.}
In this case we may take the elements of $M$ to be of the form
$j(f)(j``g(\ka))$, so in $V[G_\ka]$ we consider all $x\in\Power_\ka(g(\ka))$;
again, because of our definition of a $g(\ka)$-supercompact coding point,
the argument goes through without difficulty.  To show that
$M[G^*]$ is closed under taking $g(\ka)$-tuples, note that for any 
$g(\ka)$-tuple from $M[G^*]$ in $V[G]$, we may consider an
$g(\ka)$-tuple $t$ of names for its elements in $V[G]$, where the names are in
$M$.
All $g(\ka)$-tuples in $V[G]$ of elements of $V$ are in $V[G_{g(\ka)}]$
by closedness of the tail of the iteration, and so since $S_{g(\ka)}$ is
trivial beyond some bound below $\ka$, there is a nice name $\dot t$ for $t$
with only $g(\ka)$ elements.  
Therefore, since $M$ is closed with respect to taking $g(\ka)$-tuples in $V$,
$\dot t\in M$, and since $G^*_{\ka}$ is the same as $G_{\ka}$,
$t\in M[G^*_{j(\ka)}]$.  
But then the original $g(\ka)$-tuple of elements from $M[G^*]$ 
is in $M[G^*]$, as desired. \\

\noindent{\bf $\eta$-Extendible Cardinals.}
On the domain and range of $j$, the nontrivial part of $S$ is bounded below
$\ka$, so this is trivial.\\

\noindent{\bf $m$-Huge Cardinals.}
This is much like the $g(\ka)$-supercompact case.
The elements of $M$ may be taken to be of the form 
$j(f)(j``(j^m(\ka)))$ where the domain of $f$ is $\Power(j^m(\ka))$,
so the assumed $2^{j^m(\ka)}$-closure of $S^\ka$ 
is what we need to construct the
$M$-generic for the lifting.
Closure of $M[G^*]$ with respect to taking $j^m(\ka)$-tuples is
exactly as in the $g(\ka)$-supercompact case.\\

\noindent This completes the verification.
\hfill\qed\end{pf}

One may wonder if the restriction on which $\phi$-cardinals are preserved
(that is, only those that are not limits of $\phi$-cardinals) can be lifted.
However, some kind of restriction like this is necessary for our technique.
We are using the fact that the set of coding points is bounded below every
cardinal that we lift.  Hence, by Fodor's theorem, we cannot hope to lift
all $\phi$-cardinals from a universe where they form a stationary set in 
$\ORD$.

On the other hand, we can always mollify this problem by restricting
it to smaller and smaller classes of large cardinals.
By thinning out the class of coding point while keeping it unbounded 
in $\ORD$, the above arguments will still go through at all of the cardinals
that were previously preserved, but with new cardinals added to the
list of large cardinals that stay large.  For example, if there
are boundedly many measurable limits of measurables, and
we thin out the measurable-coding points to
only use those ``directly after'' a measurable limit of
measurables, until there are no more, then
we will still preserve all measurable cardinals that are not limits of 
measurables, but we will also preserve those measurable limits of 
measurables that are not limits of measurable limits of measurables.
If there is a proper class of measurable limits of measurables,
the ``thinned out'' class of coding points is even easier to describe:
it is simply the set of $\phi$-coding points for 
$\phi\equiv$``$\ka$ is a measurable limit of measurables''.
Indeed, this can be done not just for $\phi$ limits of $\phi$-cardinals,
but for any proper class sequence of cardinals at whose limits we don't mind
preservation failing.  So for example, we may deduce the following.

\begin{thm}
Suppose there is a proper class of $\phi_0$-cardinals, and
let $\de$ be an arbitrary ordinal.
Then a definable well-order of the universe 
may be forced while preserving all 
measurable,
$\eta$-strong for $\eta<\de$,
Woodin,
$n$-superstrong for $n\in\omega+1$,
hyperstrong,
$\ka^{+\eta}$-supercompact for $\eta<\de$,
$\eta$-extendible for $\eta<\de$, and
$m$-huge for $m\in\omega$ cardinals that are not limits of 
$\phi_0$-cardinals.
\end{thm}
\begin{pf}
We take $\phi$ to be the (size $|\de|$) disjunction of all of the
stated large cardinal properties for the sake of defining
$\phi$-bounds, but for defining the coding points we 
take $B$ to only contain those
minimal $\phi$-bounds that are minimal above a 
$\phi_0$-cardinal or are a limit of such
(and as before, take the block of coding points starting at such minimal
$\phi$-bounds to have length the next inaccessible).
The arguments from Theorem~\ref{DWOClassLCs} for each individual 
case will all go through unaffected, as long as the cardinal in question
is not a limit of $\phi_0$-cardinals.
\hfill\qed\end{pf}

Of course it makes sense to choose a very strong large cardinal property
as $\phi_0$ in this theorem,
for the simple reason that generally, stronger large cardinals \emph{are}
limits of weaker large cardinals.  
The choice can be calibrated to the tastes of the reader ---
if the assumption
of, say, a proper class of $\omega$-superstrong cardinals 
seems unpalatably strong, one
can use some other large cardinal property as $\phi_0$
and the theorem will remain true, albeit vacuous in some cases.

\ack

The author would like to thank Sy Friedman for his guidance throughout 
the course of this research,
and James Cummings and Heike Mildenberger
for many corrections and helpful comments.
The article was written while supported as a doctoral student at the
Kurt G\"odel Research Center for Mathematical Logic of the 
University of Vienna on Austrian Science Fund (FWF) project
P~16790-N04.

\bibliographystyle{elsart-num}
\bibliography{LCDWO}

\end{document}